\documentstyle[twoside]{article} 
\newtheorem{thm}{Theorem}
\newtheorem{defn}[thm]{Definition}
\newtheorem{lemma}[thm]{Lemma}
\newtheorem{prop}[thm]{Proposition}
\newtheorem{cor}[thm]{Corollary}

\newcommand{\corn}{\multicolumn{1}{r|}{1}}
\newcommand{\ner}[1]{\cline{#1-#1}}

\newcommand{\SB}[1]{\mbox{{\bf \scriptsize #1}}}

\newcommand{\comp}{\mbox{comp}}
\newcommand{\inv}{\mbox{inv}}
\newcommand{\nin}{\mbox{nin}}
\newcommand{\Roots}[1]{\mbox{Roots}(#1)}

\newcommand{\ull}{\underline{\ll}}
\newcommand{\ugg}{\underline{\gg}}
\newcommand{\fig}[2]{\begin{table}[htbp]
	\caption{#1} 
	\begin{center}\fbox{ \scriptsize \begin{tabular}#2 \end{tabular}}
	\end{center}
\end{table}}

\begin{document} 

\begin{titlepage} 
\title{Getting Results with Negative Thinking}
\author{D. Loeb\\
LaBRI, URA CNRS 1304\\
Universit\'{e} de Bordeaux I\\
33405 Talence, France \and 
E. Damiani and O. D'Antona\\
Dipartimento di Scienze dell'Informazione\\
Universit\`a di Milano\\
Milano, Italy}
\end{titlepage} 
\maketitle

\begin{abstract}
Given a universe of discourse $U$, a {\em multiset} can be thought of as a
function $M$ from $U$ to the natural numbers ${\bf N}$. In this paper, we
define a {\em hybrid set} to be any function from the 
universe $U$ to the integers ${\bf Z}$. 
These sets are called hybrid since they contain elements with either a
positive or negative multiplicity.

Our goal is to use these hybrid sets {\em as if} they were multisets in order
to adequately generalize certain combinatorial facts which are true
classically only for 
nonnegative integers. However, the definition above does not tell us much about
these hybrid sets; we must define operations on them which provide us with the
needed combinatorial structure.

In section \ref{hss}, we will define an analog of a set which can contain
either a positive or 
negative number of elements.
We will allow sums to be calculated over an
arbitrary hybrid set and in particular over ``improper'' intervals.
This will lead us in section \ref{TCSF} to the calculation of the 
{\em complementary symmetric
functions} $\comp_{n}(H)$ whose argument is a hybrid set of variables
which is at the same time a generalization of elementary and complete
symmetric functions.

In particular, as we will see in section \ref{ccs}, the complementary
symmetric function generalizes the following two classic results concerning
sequences of polynomials with {\em persistant roots}, that is to say, having
the form $p_{n}(x)=(x-a_{1})(x-a_{2})\cdots (x-a_{n})$.
\begin{eqnarray}
p_{n}(x) &=& \sum _{k=0}^{\infty} e_{k}(-a_{1},\ldots,-a_{n}) x^{n-k}
\label{hs1}\\
x^{n} &=& \sum _{k=0}^{\infty } h_{k}(a_{1}, \ldots, a_{n-k+1})
p_{n-k}(x). \label{hs2}
\end{eqnarray}
The complementary symmetric function not only allows us to express the
connection constants between two arbitrary sequences of polynomials
with persistant roots, but to also expand arbitrary rational functions
in terms of ``extended'' sequences with persistant roots and
persistant poles.

Many sequences of combinatorial interest $c_{nk}$ can be defined in
terms of connection constants. The complementary symmetric function
therefore allows us to simultaneously study them: extending them to
the case where $n$ and $k$ might be negative, demonstrating their
recursions, characterizing them with simple formulas, and providing
general inversion formulas. Along these lines, in section \ref{bcs}, we
calculate  generalized binomial coefficients and give them
a combinatorial interpretation in
terms of a 
partial order on the hybrid sets. 
Similarly, we generalize the Gaussian coefficients, and
explain their combinatorial significance. Finally, in section
\ref{sns}, we study a generalization of the $p,q$-Stirling numbers of
the first and second kind. 
\end{abstract}

\section{Hybrid Sets}\label{hss}

\subsection{Notation}
One usually thinks of a multiset $M$ as a function from some
universe $U$ to the natural numbers ${\bf N}$. $M(x)$ is
called the {\em multiplicity} of $x$ in the multiset $M$. 
Sets then are a special case of multisets in which all $x$
have multiplicity either 1 or 0.

In this paper, we generalize these concepts to
new ``sets'' and new ``multisets'' which may have negative integers as
multiplicities as well as nonnegative integers.

\begin{defn}[Hybrid Sets]\label{H-set}\label{size}
Given a universe $U$, any function $f:U\rightarrow {\bf Z}$
is called a {\em hybrid set}. 

The value of $f(x)$ is said to be the {\em
multiplicity} of the element $x$. If $f(x)\neq 0$ we say $x$ is a member of $f$
and write $x\in f$; otherwise, we write $x\not\in f$. If the support
of $f$ is finite, define the {\em number of 
elements} $\# f$ to be the sum 
$\sum_{x\in U}f(x)$. $f$ is said to be an $\# f$-element
hybrid set.

Given two hybrid sets $f$ and $g$, their {\em sum} or {\em union} is
written as 
the hybrid set $f+g$ for which the multiplicity of $x$ is the sum of
the multiplicities in $f$ and $g$: $(f+g)(x) = f(x) + g(x)$.
Similarly, the hybrid set {\em difference} of $f$ and $g$ is denoted
$f-g$ and is given by $(f-g)(x) = f(x) - g(x)$.
\end{defn}

We will denote hybrid sets by using the usual set braces $\{ \}$ and inserting
a bar in the middle $\{| \}$. Elements occuring with positive multiplicity are
written on the left of the bar, and elements occuring with negative
multiplicity are written on the right. Order is completely irrelevant.

For example, if $f=\{a,b,c,b|d,e,e \}$ then $f(a)=1$, $f(b)=2$, $f(c)=1$,
$f(d)=-1$, $f(e)=-2$, and $f(x)=0$ for $x\neq a,b,c,d,e$. 

So far, we haven't accomplished much. These hybrid sets are equivalent
to elements of the free group over the universe in question, so one
might wonder if there is anything new here. However, we are planning
on treating these hybrid sets as if they were multisets. That is to
say, we will consider the subsets of hybrid sets; we will perform
sums over hybrid sets; and we will enumerate certain classes of
combinatorial objects with them. We begin by considering which hybrid
sets correspond to the classical notion of a set.  

\begin{defn}[New Sets]
A {\em positive} or {\em classical set} is a hybrid set taking only values $0$
and $1$. A {\em negative set} is a hybrid set taking only values $0$ and $-1$.
A {\em new set} is either a positive or a negative set.
\end{defn}

For example, the set $S=\{a,b,c \}$ can be written as the classical or positive
set $S = \{b,a,c| \}$ while
$-S=\{|a,b,c \}$ is a negative set. Both $S$ and $-S$ are new sets. However,
$\{a|b \}$ is not a new set.  

The {\em empty set} $\emptyset=\{| \}$ is the unique hybrid set for which  all
elements have multiplicity zero. It is thus simultaneously a positive set and a
negative set.

\subsection{Sums and Products over a Hybrid Set}

One often wishes to sum an expression over all the elements
of a {\em set}. For example, we usually write 
$$ \phi(n) = \sum_{i\in \mbox{fact}(n)}i $$
for the Euler phi function where $\mbox{fact}(n)$ is the set of
positive factors of an integer $n$.

When dealing with {\em multisets}, it is often convenient to
include multiplicities in the sum or product. One defines
$ \sum_{x\in M} F(x)$ and $\prod _{x\in M} F(x)$ to be $\sum_{x\in
U}M(x)F(x)$ and $\prod _{x\in U}F(x)^{M(x)}$ respectively. For
example, we can expand a monic polynomial $p(x)$ by
$$ p(x) = \prod_{\lambda \in \Roots{p}} (x-\lambda ) $$
where $\Roots{p}$ is the multiset of roots of the polynomial $p$.

This notation can be immediately generalized in the obvious way to
{\em hybrid sets}. 
\begin{defn}{\bf (Sums and Products over Hybrid Sets)}
For any function $F$ and hybrid set $h$ we write
$\displaystyle \sum_{x\in h}F(x)$ for $\displaystyle \sum_{x\in
U}h(x)F(x)$, and we write 
$\displaystyle \prod_{x\in h}F(x)$ for $\displaystyle \prod_{x\in
U}F(x)^{h(x)}$. 
\end{defn}

For example, consider the following extension of the definition of
$\Roots{p}$.
\begin{defn}{\bf (Hybrid Set of Roots)}
Let $f(x)$ be a rational function. Then $f(x)$ can be written as
$$ f(x)=c\frac{(x-a_{1})(x-a_{2})\cdots
(x-a_{i})}{(x-b_{1})(x-b_{2})\cdots (x-b_{j})}.$$
We then call
$$ \Roots{f} = \{ a_{1},a_{2},\ldots,a_{i} | b_{1},b_{2}, \ldots,b_{j}
\} $$
the {\em Hybrid Set of Roots} of $f(x)$
\end{defn}
For $f(x)$ a monic rational function, we thus have
$$ f(x) = \prod _{\lambda \in \Roots{f}} (x-\lambda ). $$

\subsection{Sums and Products with Limits}

Often, sets are written using ellipsis for missing elements. For
example, $\{a_{1},\ldots,a_{5} \}$ stands for the set
$\{a_{1},a_{2},a_{3},a_{4},a_{5} \}$. In the case of hybrid sets, this
notation turns out to be rather powerful.
\begin{defn}{\bf (Ellipsis)}
Let $(a_{n})_{n\in \SB{Z}}$ be a sequence of constants indexed by an
integer $n$. Their {\em ellipsis} is defined by the recursion:
\begin{eqnarray*}
\{a_{1}..a_{0} \} &=& \emptyset = \{| \}\\
\{a_{1}..a_{n} \} &=& \{a_{1}..a_{n-1} \} + \{a_{n} | \}.
\end{eqnarray*}
\end{defn}
It therefore follows that $\{ a_{i}..a_{j} \}$ is
explicitly given 
by: 
$$ \{a_{i}..a_{j} \}=\left\{ \begin{array}{ll}
\{a_{i},a_{i+1},\ldots ,a_{j-1} ,a_{j}|\}&\mbox{if $i \leq j$,}\\[0.1in]
\emptyset &\mbox{if $i=j+1$, and }\\[0.1in]
\{|a_{j+1},a_{j+2},\ldots ,a_{i-1} \}&\mbox{if $i\geq j+2$;}
\end{array} \right. $$
In the first case, we are reduced to the standard usage of ellipsis.
However, in the other cases, the interval being considered is
``improper'' classically, but the meaning is clear here. 
Reference is
being made to a certain $j-i+1$-element negative set. Our ellipsis are
written with only two dots to distinguish them from normal ellipsis.

One often uses the notation $\displaystyle \sum_{n=1}^{i}A(n)$ where
$n$ is a nonnegative integer to denote the sum of $A(n)$ over
the set $\{1,2,\ldots ,n \}$. Using the ellipsis notation,
we can define the sum or product of
a  quantity from an arbitrary integer to another:  
\begin{defn}{\bf (Sums and Products with Limits)}
Let $(A(n))_{n\in \SB{Z}}$ be a sequence of constants indexed by an
integer $n$, and let $i$ and $j$ be arbitrary integers. The sum and
product of $A(n)$ from $i$ to $j$ is given by
\begin{eqnarray*}
 \sum _{n=i}^{j} A(n) &=& \sum _{n\in \{i..j \}} A(n) \\
 \prod _{n=i}^{j} A(n)&=& \prod _{n\in \{i..j \}} A(n) .
\end{eqnarray*}
\end{defn}
Just as with integrals, the sum is shift invariant, the sum over an empty
interval is zero, and the sum over a positive expression over an ``improper''
interval $j<i$ is negative. 

For example, $\sum _{i=1}^{4}i=1+2+3+4=10$ since $\{1..4 \}=\{1,2,3,4|
\}$, and $\sum_{i=1}^{-4}=-(-3)-(-2)-(-1)-0=6$ since $\{1..-4 \} = \{|
0,-1,-2,-3 \}$. 
In general, $\sum _{i=1}^{n}i={n(n+1)}/{2}$ for all $n$ regardless of
sign. In fact,

\begin{prop}
Let $p(x)$ be any polynomial, then there is a polynomial $q(x)$ such that
$q(n)=\sum_{i=1}^{n}p(i)$ for all $n$ regardless of sign.$\Box $ 
\end{prop}

{\em Proof:} Classical proof by induction can be applied here {\em mutatis
mutandis}.$\Box$ 

Finally, note the following reciprocity formula useful in calculations
with improper intervals:
\begin{prop}
Let $(a_{i})_{i\in \SB{Z}}$ be a sequence of constants indexed by an
integer $i$, then for all integers $n$
$$ \{a_{1}..a_{n} \} = -\{b_{1}..b_{-n} \}$$
where $b_{i}=a_{n+i}.\Box $
\end{prop}

\section{The Complementary Symmetric Function} \label{TCSF}
\subsection{Introduction}

The complementary symmetric function is so called because it {\em
complements} or extends the {\em compl}ete and el{\em ementary}
symmetric functions.
Recall that $e_{n}(V)$ and $h_{n}(V)$, the elementary and complete
symmetric functions of degree 
$n$ over the set of variables $V$, are given by the following
generating functions:
\begin{eqnarray*}
\sum _{n\geq 0}h_{n}(V)t^{n} &=& \prod _{x\in V}(1-xt)^{-1}\\
\sum _{n\geq 0}e_{n}(V)t^{n} &=& \prod _{x\in V}(1+xt).
\end{eqnarray*}
The {\em complementary symmetric function} combines these two
definitions by using our notion of the product over a hybrid set.
\begin{defn}
Let $V$ be a hybrid set of variables. Then the generating function for
the complete symmetric function of degree $n$ over these variables is
given by the product:
$$ \sum _{n\geq 0} \comp_{n}(V)t^{n} = \prod _{x\in V}(1-xt). $$
\end{defn}
In other words, we have
$$\sum _{n\geq 0} \comp_{n}(a_{1},a_{2},\ldots|b_{1},b_{2},\ldots)
t^{n} =  \frac{\displaystyle \prod _{i=1}^{\infty }(1-a_{i}t)} 
{\displaystyle \prod _{i=1}^{\infty}(1-b_{i}t)}
= \frac{\displaystyle \sum _{n=0}^{\infty } e_{n}(-a_{1},-a_{2},\ldots)t^{n}}
{\displaystyle \sum _{n=0}^{\infty } h_{n}(a_{1},a_{2},\ldots)t^{n}}.$$
In the case where $V$ is a negative or positive set, it follows that
\begin{eqnarray}
e_{n}(a,b,c,\ldots)=\comp_{n} (-a,-b,-c,\ldots|) \label{p2a}\\
h_{n}(a,b,c,\ldots)=\comp_{n} (|a,b,c,\ldots). \label{p2b}
\end{eqnarray}
It is now immediate that the transformation
$$ *: e_{n}\mapsto h_{n} $$
used in the theory of symmetric functions is an involution. $*$ acts
on the set of symmetric variables, inverting the sign and the
multiplicity of each variable.

The next few results are based on the following convolution identity.
\begin{lemma}\label{l2}
Suppose $V$ and $W$ are two (hybrid) sets of variables. Then
$$ \comp _{n}(V+W) = \sum _{k=0}^{n} \comp _{k}(V) \comp _{n-k}(W). $$
\end{lemma}

{\em Proof:} Using the definition of $V+W$ (Definition~\ref{H-set}),
we can proceed as follows
\begin{eqnarray*}
\sum _{n=0}^{\infty } \comp _{n}(V+W) t^{n} &=& \prod _{x\in V+W}(1-xt)\\
 &=&\left( \prod _{x\in V}(1-xt) \right)\left( \prod _{x\in W}(1-xt)\right) \\
&=& \left(\sum _{n=0}^{\infty } \comp _{n}(V) t^{n}  \right)
\left(\sum _{n=0}^{\infty } \comp _{n}(V) t^{n}  \right)\\
&=& \sum _{n=0}^{\infty }\left(\sum _{k=0}^{n}\comp _{k}(V) \comp
_{n-k}(W)\right) t^{n},
\end{eqnarray*}
and identify coefficients of $t^{n}.\Box$

We have the following explicit characterization of the complementary
symmetric function.
\begin{thm}\label{t3.1}
The complementary symmetric function $\comp
_{n}(a_{1},a_{2},\ldots|b_{1},b_{2},\ldots)$ is the sum of all
monomials in the $a$'s and the $b$'s in which the $a$'s may appear
with no exponent two or greater, and are associated with a minus sign.
\end{thm}

For example, 
\begin{eqnarray*}
\comp _{2}(a,b,c|x,y,z) &=& ab+ac+bc\\
&	&	-ax-ay-az-bx-by-bz-cx-cy-cz\\
&	&	+x^{2}+xy+xz+y^{2}+yz+z^{2}.
\end{eqnarray*}

{\em Proof of Theorem \ref{t3.1}:} First, by Lemma \ref{l2},
$$ \comp _{n}(a_{1},a_{2},\ldots|b_{1},b_{2},\ldots) = \sum _{k=0}^{n}
\comp _{k}(a_{1},a_{2},\ldots|) \comp _{n-k}(|b_{1},b_{2},\ldots). $$
Then by equations \ref{p2a} and \ref{p2b}, we have in turn
$$ \comp _{n}(a_{1},a_{2},\ldots|b_{1},b_{2},\ldots) = \sum _{k=0}^{n}
e_{k}(-a_{1},-a_{2},\ldots) h_{n-k}(b_{1},b_{2},\ldots) $$
where as we recall the elementary symmetric function is the sum of
monomials with no repeated variables, and the complete symmetric
function is the sum of monomials possibly with repeated
variables.$\Box $

The complementary symmetric function satisfies an interesting
recursion:
\begin{thm}\label{t5}
For all $n,m,k$ positive integers, we have
\begin{eqnarray*}
 \comp _{n}(a_{1},\ldots,a_{m}|b_{1},\ldots,b_{k}) &=&
 \comp _{n}(a_{1},\ldots,a_{m-1}|b_{1},\ldots,b_{k-1}) \\
&& +(b_{k}-a_{m})
\comp_{n-1}(a_{1},\ldots,a_{m-1}|b_{1},\ldots,b_{k}).\Box 
\end{eqnarray*}
\end{thm}

{\em Proof:} Either apply Lemma \ref{l2} first with $W=\{a_{m}|b_{k}
\}$ and then again with $W=\{a_{m}| \}$, or reason with Equations
\ref{p2a} and \ref{p2b}.$\Box $

Note that Theorem \ref{t5} holds equally well for $m,k\leq 0$ provided
that $\comp _{n}(a_{1},\ldots,a_{m}|b_{1},\ldots,b_{k})$ is read as the
complementary symmetric function
$\comp _{n}(\{a_{1}..a_{m} \} - \{b_{1}..b_{k} \})$ whose argument is
the difference of two hybrid sets each given by ellipsis.

The following expansion can be found in [3].

\begin{thm}
The complementary symmetric function is given by the following sum
$$ \comp _{n}(a_{1},a_{2},\ldots|b_{1},b_{2},\ldots) = \sum
_{0<\alpha_{1}<\alpha_{2}<\cdots <\alpha_{n}} \prod _{i=1}^{n}
(b_{\alpha_{i}+1-i}-a_{\alpha_{i}}) $$
over strictly increasing sequences of positive integers.
\end{thm}

{\em Proof:} Apply Theorem \ref{t5} until all complementary symmetric
function disappear either by reducing their degree to
zero, or by eliminating all their variables.$\Box $

\subsection{Connection Constants}\label{ccs}

In order to state the main result of this section, we must first
introduce the concept of persistant roots. 
\begin{defn}
Let $(b_{n})_{n\in \SB{Z}}$ be a sequence of constants indexed by an
integer, and let 
$$ q_{n}(x) = \prod _{i=1}^{n} (x-b_{i}) $$
for all $n$ positive or negative. In other words,
\begin{eqnarray*}
\vdots &&  \vdots\\
q_{2}(x) &=& (x-b_{1})(x-b_{2})\\
q_{1}(x) &=& (x-b_{1})\\
q_{0}(x) &=& 1\\
q_{-1}(x) &=& 1/(x-b_{0})\\
q_{-2}(x) &=& 1/(x-b_{0})(x-b_{-1})\\
\vdots && \vdots
\end{eqnarray*}
Then $(q_{n}(x))_{n\in \SB{Z}}$ is called the {\em sequence of
rational functions with persistant roots} $(b_{n})_{n\in \SB{Z}}$.
\end{defn}

We can now state the main result concerning the complementary
symmetric function.
\begin{thm}\label{t9}
Let $f(x)$ be a monic rational function of degree $n$, and
$(b_{k})_{k\leq n}$ a sequence of constants indexed by an integer.
Then $f(x)$ can be expanded in terms of the sequence of rational
functions $(q_{k}(x))_{k\leq n}$ with persistant roots $(b_{k})_{k\leq
n}$ as follows:
\begin{equation}\label{e*}
f(x)=\sum _{k=0}^{\infty } \comp _{k}(\Roots{f}-\{b_{1}..b_{n-k+1}
\}) q_{n-k}(x).
\end{equation}
\end{thm}

Notes to {\em Theorem \ref{t9}:} \begin{enumerate}
\item The degree of a rational function $f(x)=p(x)/q(x)$ is defined to
be the degree of the numerator $p(x)$ minus the degree of the
denominator $q(x)$.
\item For $0\leq k\leq n$, the coefficient of $q_{k}(x)$ does not
depend on $b_{k+2},b_{k+3},\ldots,b_{n}$ even though initial
consideration of the problem might seem to indicate that it should. 
\item To interpret this theorem when $k>n$, recall the definition of
improper intervals from section 3.
\item The coefficient of $q_{n}(x)$ does not depend on the
non-existant ``$b_{n+1},$'' since $\comp _{0}(V)$ is identically one.
\item The case in which $f(x)$ is a polynomial (see Corollary
\ref{c11}) is much simpler and is a direct consequence of equations
\ref{hs1}, \ref{hs2}, \ref{p2a}, \ref{p2b}, and Lemma \ref{l2}.
\end{enumerate}

{\em Proof of Theorem \ref{t9}:} We proceed by induction on the number<
of roots and poles 
which $f(x)$ possesses. The base case is thus $f(x)=1$. In this case,
the right hand side of equation \ref{e*} is
$$ \sum _{k=0}^{\infty }\comp _{k}(-\{b_{1}..b_{1-k}\})q_{-k}(x). $$
For $k=0$, the summand is one, since the complementary symmetric
function of degree zero is identically one and $q_{0}(x)=1$. The
remaining terms are all of the form
$$
\frac{e_{k}(-b_{0},-b_{-1},\ldots,-b_{2-k})}{(x-b_{0})(x-b_{-1})\cdots
(x-b_{1-k})} $$ 
which is zero for lack of sufficiently many arguments to the
elementary symmetric function. Thus, the sum is one as required.

Next, towards induction assume Eq. \ref{e*} holds for $g(x)$. Now,
consider $f(x)=(x-a)g(x)$.
\begin{eqnarray*}
f(x) &=& (x-a) \sum _{k=0}^{\infty } \comp
_{k}(\Roots{g}-\{b_{1}..b_{n-k} \}) q_{n-k}(x) \\
&=&  \sum _{k=0}^{\infty } \comp
_{k}(\Roots{g}-\{b_{1}..b_{n-k} \}) [(x-b_{n-k})+(b_{n-k}-a)]
q_{n-k-1}(x) \\ 
&=& q_{n}(x) + \sum_{k=0}^{\infty } [ \begin{array}[t]{ll}
     \comp _{k+1}(\Roots{g}-\{b_{1}..b_{n-k-1} \})\\
     + (b_{n-k}-a)\comp _{k}(\Roots{g}-\{b_{1}..b_{n-k} \})&]q_{n-k-1}(x)
     \end{array}\\
&=& q_{n}(x) + \sum _{k=0}^{\infty } \comp_{k}
(\Roots{f}-\{b_{1}..b_{n-k} \}) q_{n-k-1}(x) \\ 
&=& \sum _{j=0}^{\infty } \comp_{k}
(\Roots{f}-\{b_{1}..b_{n-k-1} \}) q_{n-k}(x)
\end{eqnarray*}
where the fourth equality results from the recursion theorem
(\ref{t5}). 

The proof for $f(x)=g(x)/(x-a)$ can be found by reading the above
sequence of equalities in reverse order.$\Box $

Theorem \ref{t9} gives rise to the following inversion formula.
\begin{cor}\label{c10}
Let $(a_{n})_{n\in \SB{Z}},$  $(b_{n})_{n\in \SB{Z}},$
$(c_{n})_{n\in \SB{Z}},$ and $(d_{n})_{n\in \SB{Z}}$  be sequences
of real numbers. Then
$$ c_{n} = \sum _{k=0}^{\infty } \comp _{k}(\{a_{1}..a_{n} \}-
\{ b_{1}..b_{n-k+1}\}) d_{n-k}$$ 
for all $n$, if and only if
$$ d_{n} = \sum _{k=0}^{\infty } \comp _{k}(\{b_{1}..b_{n} \}-
\{\a_{1}..a_{n-k+1}\}) c_{n-k}$$ 
for all $n.\Box$
\end{cor}

Theorem \ref{t9} is particularly easy to state in the case of
polynomials. 
\begin{cor}\label{c11}
Given constants $(a_{i})_{i=1}^{\infty }$ and $(b_{i})_{i=1}^{\infty
}$, we have 
$$ (x-a_{1})\cdots (x-a_{n}) = \sum _{k=0}^{n} \comp
_{k}(a_{1},\ldots,a_{n}|b_{1},\ldots,b_{k+1}) (x-b_{1})\cdots
(x-b_{k}).\Box $$
\end{cor}

As special cases, we have Eqs. \ref{hs1} and \ref{hs2}. Other special
cases of particularly interesting connection constants along with
their combinatorial interpretation in the section below and throughout
the rest of this paper.

\subsection{Example of Connection Constants}

We begin with an example of the expansion of an inverse formal power
series since this lies outside of the usual theory. Let us expand
$1/x$ in terms of the sequence of rational functions $q_{n}(x)$ with
persistant poles $b_{0}, b_{-1}, b_{-2}, \ldots$; that is, expand
$1/x$ with respect to
\begin{eqnarray*}
q_{-1}(x)&=& 1/(x-b_{0}),\\
q_{-2}(x)&=& 1/(x-b_{0})(x-b_{-1}),\\
q_{-3}(x)&=& 1/(x-b_{0})(x-b_{-1})(x-b_{-2})\mbox{, etc.}
\end{eqnarray*}

By theorem \ref{t9}, we can write
\begin{eqnarray*}
\frac{1}{x} &=& \sum _{k= 0}^{\infty } \comp
_{k}(|0,b_{0},b_{-1},\ldots,b_{1-k}) q_{-1-k}(x)\\ 
&=& \sum _{k= 0}^{\infty } e_{k}(-b_{0},-b_{-1},\ldots,b_{1-k})
q_{-1-k}(x)\\ 
&=& \sum _{k= 0}^{\infty } (-1)^{k} b_{0}b_{-1}\cdots b_{1-k}
q_{-1-k}(x).
\end{eqnarray*}
So in particular,
$$ \frac{1}{x} = \sum _{n=1}^{\infty } \frac{1}{(x+1)^{n}} $$
for $x>0$ as the convergent sum of a geometric series. Similarly,
\begin{equation}\label{telescope}
\frac{1}{x} = \sum _{n=1}^{\infty } \frac{(k-1)!} {(x+1)
(x+2)\cdots (x+k)}.
\end{equation}
Using the classical method of telescoping series, equation
\ref{telescope} can be validated for $x$ a positive integer. However,
to demonstrate this 
identity in general for $x>0$, the preceeding theory of connection
constants is necessary.

Other examples of connection constants can be found in the following
sections.

\section{Binomial Coefficients} \label{bcs}
\subsection{Connection Constants}
Another classical example of connection constants is that of the
binomial coefficient ${n\choose k}$ which is given by the coefficients
of the polynomial $(x+1)^{k}$, or more generally the Gaussian
coefficients ${n\choose 
k}_{q}$ which is given by the coefficients of
the polynomial $(x-1)(x-q)(x-q^{2})\cdots (x-q^{n-1})$ (see Figure
\ref{fig:Gauss1} and \ref{fig:Gauss2}). These
constants then have 
well known combinatorial interpretations. ${n\choose k}$ is the number
of $k$-element subsets of a given $n$-element set whereas the
coefficient of $q^{t}$ in ${n\choose k}_{q}$ 
counts the number of partitions $\lambda \vdash t$ of length $k$ with
all parts distinct and less than or equal to $n$ where a partition
of $t$ is defined as a 
nonincreasing sequence of nonnegative integers summing to
$t$.\footnote{Our definition of the Gaussian coefficient differs from with
certain other definitions (for example, see \cite{Leroux}) by a factor
of $q^{(n-k}(n-k-1)/2$. Under this alternate definition, the
coefficients of ${n\choose k }_{q}$ then count the number of
partitions $\lambda \vdash t$ of length $k$ with 
all parts less than or equal to $n-k$. Moreover for $q$ prime, ${n
\choose k}_{q}$ counts the
number of $k$-dimensional subvector spaces in an $n$-dimension vector
spaces over a $q$-element field.}

Now, that we are capable of calculating these connection constants via
Theorem \ref{t9}, we can equally well interpret them for $n$ negative
or positive. We are then calculating the coefficient of $x^{k}$ in an
inverse Laurent series.

By Theorem \ref{t5}, we have the recurrences
\begin{eqnarray}
{n\choose k }_{q} &=& {n-1\choose k-1 }_{q} + q^{n-1}{n-1\choose k }_{q}
\label{rec1}\\
{n\choose k } &=& {n-1\choose k-1 } + {n-1\choose k }. \label{rec2}
\end{eqnarray}

Moreover, by theorem \ref{t9}, we have
\begin{eqnarray}
{n\choose k}_{q} &=& \comp _{n-k}  (-1, -q, ..,-q^{n-1}) \nonumber \\
&=& \left\{\begin{array}{ll}
    e_{n-k}(1,q,\ldots,q_{n-1})&	
\mbox{for $0\leq k\leq n$, and}\\[0.12in] 
    (-1)^{n-k} h_{n-k}(q^{-1},q^{-2},\ldots,q^{n})&	\mbox{for
    $k\leq n<0$}.
\end{array} \right.  \label{heq} \\
&=& \frac{\displaystyle \prod _{i=1}^{k}q^{n+1-i}-1}
{(q^{k}-1)(q^{k-1}-1)\cdots (q-1)}. \nonumber
\end{eqnarray}
And in particular for $q=1$,
\begin{eqnarray*}
{n\choose k} &=& \comp _{n-k}  (\underbrace{-1, -1,
..,-1}_{n\mbox{\tiny{} terms}})\\[0.1in]
&=& \left\{\begin{array}{ll}
    e_{n-k}(\underbrace{1,1,\ldots,1}_{n\mbox{\tiny{} terms}})&
\mbox{for $0\leq k\leq n$, and}\\[0.45in] 
    (-1)^{n-k} h_{n-k}(\underbrace{1,1,\ldots,1}_{-n\mbox{\tiny{} terms}})&
\mbox{for $k\leq n<0$}
\end{array} \right.\\[0.1in]
&=& \frac{n(n-1)(n-2)\cdots (n-k+1)}{k!}.
\end{eqnarray*}
Setting all the variables of the elementary or complete symmetric
function to one   gives the
number of terms. In the case of the elementary symmetric function,
this is the number of ways of choosing $k$ variables among $n$ without
replacement. Thus, we have rediscovered the classical combinatorial
interpretation of ${n\choose k}$ for $n$ positive. For the complete
symmetric function, the number of terms is the number of ways of
choosing $k$ variables among $n$ with replacement. In other words, up
to sign, ${-n\choose k}$ is  the number of $k$-element multisets based
on a given $n$-element set. Hence, we have combinatorial
interpretations of ${n \choose k}$ for $n$ positive or negative.

\subsection{Roman Coefficients}
What is the connection between these two interpretations? And can we
permit $k$ to be negative as well?

We can for example, adopt the following definition of ${n\choose k}$
(consistant with the definition above) for all values of $n$ and $k$,
positive or negative (see Figure~\ref{fig:Bin}).
\begin{equation}\label{romco}
 {n \choose k}=\lim_{\epsilon \rightarrow 0}\frac{\Gamma
(n+1+\epsilon)}{\Gamma 
(k+1+\epsilon ) \Gamma(n-k+1+\epsilon )}. 
\end{equation}

\fig{Binomial Coefficients, ${n\choose k}$}{{r|rrrr|*{7}{r}}
\label{fig:Bin}
\normalsize $n \backslash k$ &\normalsize  $-$4 &\normalsize  $-$3
&\normalsize  $-$2 &\normalsize  $-$1 &\normalsize  0 &\normalsize  1
&\normalsize  2 &\normalsize  3 &\normalsize  4 &\normalsize
\normalsize  5 &\normalsize  6\\
\hline  
\normalsize 6 & 0&0&0&0 & 1 & 6 & 15 & 20 & 15 & 6 & \corn \\[0.5mm] \ner{12}
\normalsize 5 & 0&0&0&0& 1 & 5 & 10 & 10 & 5 & \corn  &0\\[0.5mm]  \ner{11}
\normalsize 4 & 0&0&0&0& 1 & 4 & 6 & 4 & \corn  &0&0\\[0.5mm] \ner{10}
\normalsize 3 & 0&0&0&0& 1 & 3 & 3 & \corn  &0&0&0\\[0.5mm] \ner{9}
\normalsize 2 & 0&0&0&0& 1 & 2 & \corn  &0&0&0&0\\[0.5mm] \ner{8}
\normalsize 1 & 0&0&0&0& 1 & \corn  & 0&0&0&0&0\\[0.5mm] \ner{7}
\normalsize 0 & 0&0&0&0& \corn  & 0&0&0&0&0&0\\[0.5mm]
\hline 
\normalsize $-$1& $-$1 & 1 & $-$1 & 1 & 1 & $-$1 & 1 & $-$1 & 1 & $-$1
& 1\\[0.5mm] \ner{5}
\normalsize $-$2& 3 & $-$2 & \corn & 0 & 1 & $-$2 & 3 & $-$4 & 5
&$-$6&7 \\[0.5mm] \ner{4}
\normalsize $-$3& $-$3 & \corn  & 0&0& 1 & $-$3 & 6 & $-$10 &
15&$-$21&28 \\[0.5mm] \ner{3}
\normalsize $-$4& \corn  & 0&0&0& 1 &$-$4 & 10 & $-$20 & 35&$-$56&84
\\[0.5mm] \ner{2}
\normalsize $-$5& 0&0&0&0& 1 & $-$5 & 15 & $-$35 & 70&$-$126&210}

\begin{prop}[The Six Regions]\label{regions}
Let $n$ and $k$ be integers. Depending
on what region of the Cartesian plane $(n,k)$ is in, the following
formulas apply:
\begin{enumerate}
\item If $n\geq k\geq 0$, then ${n\choose k }= \frac{n!}{k!(n-k)!}.$
\item  If $k\geq 0>n$, then ${n\choose k } =(-1)^{k}{-n+k-1\choose k }.$
\item  If $0>n\geq k$, then ${n\choose k }=(-1)^{n+k}{-k-1\choose n-k }$
\item If $k>n\geq 0$, ${n\choose k }=0. $
\item If $0>k>n$, ${n\choose k }=0. $
\item If $k\geq 0>n$, ${n\choose k }=0. \Box $
\end{enumerate}
\end{prop}

Note that in regions 4--6 there is an extra factor of $\epsilon $ in the
numerator of the limit, so we are left with zero. Next, region 1 is the
classical case, so we have the usual binomial coefficients. Finally,
regions 2 and 3 are simply  related to region 1.

Most of the usual properties of binomial coefficients hold true in all six
regions. 

\begin{prop}[Complementation] \label{compl}
For all integers $n$ and $m$,
${n\choose m}={n\choose n-m }.\Box $
\end{prop}

\begin{prop}[Iteration]
For all integers $i,j,$ and $k$,
$$ {i\choose j }{j\choose k }= {i\choose k }{i-j\choose j-k }. \Box $$
\end{prop}

By equation \ref{rec2}, 

\begin{prop}[Pascal] \label{pascal}
Let $n$ and $k$ be integers not both zero, then
$$ {n\choose k}={n-1\choose k }+{n-1\choose k-1 }. \Box $$
\end{prop}

Nevertheless, ${0\choose 0 }=1$ while
${ -1\choose -1 }+{ -1\choose 0 }= 1+1=2$.

The coefficients in Regions 1 and 3 have already been seen to be
connection constants; that is to say, they are the coefficients of
$x^{n-k}$ in the inverse power series $(1+x)^{n}$. Those in regions 4--6
are zero. This leaves only 
region 3, which as we see below, can be interpretted as a set of
connection constants in terms of formal power series.

\begin{prop}\label{pr3}
For all integers $n$,
$(x+1)^{n}$ is given by
$$ \sum _{k=-\infty }^{n} {n \choose k}x^{n-k} $$
as a formal power series. 
\end{prop}

{\em Proof:} The $(n-k)$th derivative of $(1+x)^{n}$ evaluated at zero is
$(n-k)!{n\choose k}.\Box $

\subsection{Inclusion of Hybrid Sets}

These new binomial coefficients ${n\choose k}$ are always integers,
but what do they count? In region 1, they count the number of
$k$-element subsets of a  
given $n$-element set. We claim that given a suitable generalization of the
notion of a subset, this is true in general (up to sign).
There is no need for a separate interpretation of ${n \choose k}$ in
regions 2 and 3.

We define a partial order on these sets which is a generalization of the
usual ordering of classical sets and multisets by inclusion.
\begin{defn}[Subsets]
Informally, we say
$f$ is a {\em subset} of $g$ (and write $f\subseteq g$) if one can
remove elements 
one at a time from $g$ (never 
removing an element that is not a member of $g$) and thus either acheive $f$ or
have removed $f$.
\end{defn}

For example, we might start with the hybrid set $f=\{a,b,c,c|d,e \}$. We will
remove a few of its elements one at a time. Suppose we remove $b$; this leaves
$\{a,c,c|d,e \}$. Now, $b$ is no longer an element, so we can not remove it
again. Instead, we might remove $d$ leaving $\{a,c,c|d,d,e \}$. Obviously, we
can remove $d$ as many times as we choose. Finally,
suppose we remove $c$ leaving $\{a,c|d,d,e \}$. Hence, we have proven two
things. Since we were able to remove $\{b,d,c| \}$, we know that $\{b,d,c |
\}\subseteq f$. Also, since we were left with $\{a,c|d,d,e \}$, we know that
$\{a,c|d,d,e \}\subseteq f$.

\begin{prop}\label{patl} \cite{neg}
The ordering ``$\subseteq$'' defined above is a well defined partial
order.$\Box $
\end{prop}

{\em Proof:} See \cite{neg}.

Note however that the ordering $\subseteq$ does not form a lattice.
For example,  hybrid 
sets $f=\{|a,b \}$ and $g=\{a|b \}$ have lower bounds $\{a| \}$, $\{a,b| \}$,
and $\{a,b,b| \}$, but no greatest lower bound. Conversely, $f$ and $g$ possess
no upper bounds.

\begin{thm}\label{clas}
The subsets of a classical set $f$ correspond to the classical subsets of 
a classical set. The subsets of a multiset $f$ correspond to the
classical submultisets of a multiset.
\end{thm}

{\em Proof:} To construct a conventional subset of a set $S$, we
merely  remove some  elements subject to the
conditions that we only remove elements of $S$ and we don't
remove an element twice. The order of removal is not relevant.

To construct a conventional submultiset of a multiset $M$, we
merely  remove some  elements subject to the
conditions that we do not remove any element more times than 
its multiplicity in $M.\Box $

Now, we have the desired interpretation of the binomial coefficients.

\begin{thm}
Let $n$ and $k$ be arbitrary integers.
Let $f$ be an $n$-element new set. Then $|{n\choose k}|$ counts the
number of $k$-element hybrid sets which are subsets of $f.\Box $ 
\end{thm}

Let us consider this result in each of the six regions mentioned in
Proposition 
\ref{regions}. 
\begin{enumerate}
\item This is the only classical case. In this region, one might count the
number of 2-element subsets of the set $\{a,b,c,d| \}$. By Theorem~\ref{clas},
we enumerate the usual subsets: $\{a,b| \}, \{a,c| \}, \{a,d| \}, \{b,c| \},
\{b,d |\},$ and $ \{c,d| \}$ but no others.
\item In this region, one might count the number of 2-element subsets of
$f=\{|a,b,c \}$. These subsets correspond to what we can remove from $f$, since
what we would have left over after a removal would necessarily contain a
negative number of elements. We can remove any of the three elements any number
of times, so we have: $\{a,a |\}, \{b,b| \}, \{c,c |\},\{a,b| \}, \{b,c| \},$
and $\{a,c |\}$. 
\item Here, we are interested in $-5$-element subsets of $f$. Since
$f$ contains $-3$ elements, we must start with $f$ and remove 2
elements. Thus,  there is one 
subset here for each subset in the corresponding position in region 2. In this
case they are: $\{|a,a,a,b,c \},\{|a,b,b,b,c \},\{|a,b,c,c,c \}, \{|a,a,b,b,c
\},\{|a,b,b,c,c \},$ and $\{|a,a,b,c,c \}$.
\item By Theorem \ref{clas}, there are no 6-element subsets of the set
$\{a,b,c,d| \}$. Once you remove 4 elements, you can not remove any more.
\item There are no $-2$-element subset of a $-3$-element set $f$. If we remove
elements from $f$, we are left with less than $-3$ elements, and have removed a
positive number of elements. In neither case have we qualified a $-2$-element
hybrid set to be a subset of $f$.
\item Again by Theorem \ref{clas}, there are no $-2$-element subsets of the set
$\{a,b,c,d| \}$, since you are not allowed to introduce elements with
negative multiplicities.
\end{enumerate}

\section{Linear Partitions}

\fig{Gaussian Coefficients, ${n\choose k}_{q}$---Region 1}{{r|*{6}{l}}  
\label{fig:Gauss1}
\normalsize $n \backslash k$ & \normalsize
0 &\normalsize  1 &\normalsize  2 &\normalsize  3 &\normalsize  4\\
\hline  
\normalsize 5 & $q^{10}$ & $q^{6}+q^{7}+q^{8}$ &
$q^{3}+q^{4}+2q^{5}+2q^{6}$ &
$q+q^{2}+2q^{3}+2q^{4}$  &  $1+q+q^{2}$\\
	&&$+q^{9}+q^{10}$ &$+2q^{7}+q^{8}+q^{9}$ &
$+2q^{5}+q^{6}+q^{7}$ & $+q^{3}+q^{4}$\\[1mm]
\normalsize 4 & $q^{6}$ & $q^{3}+q^{4}+q^{5}+q^{6}$ &
$q+q^{2}+2q^{3}+q^{4}+q^{5}$ & $1+q+q^{2}+q^{3}$ & 1 \\[0.5mm] 
\normalsize 3 & $q^{3}$ & $q+q^{2}+q^{3}$& $1+q+q^{2}$ & 1 &0\\[0.5mm]
\normalsize 2 & $q$ & $1+q$ & 1 &0&0\\[0.5mm]
\normalsize 1 & 1 & 1 & 0&0&0\\[0.5mm]
\normalsize 0 & 1 & 0&0&0&0}
\fig{Gaussian Coefficients, ${n\choose k}_{q}$---Region 2}{{r|lllll}  
\label{fig:Gauss2}
\normalsize $n \backslash k$ &\normalsize $-5$  &\normalsize  $-$4
&\normalsize  $-$3 &\normalsize  $-$2 &\normalsize  $-$1 \\
\hline  
\normalsize $-$1& $q^{-4}$ & $-q^{-3}$ & $q^{-2}$ & $-q^{-1}$ & 1\\[0.5mm]
\normalsize $-$2& $-q^{-6}-q^{-5}-q^{-4}-q^{-3}$ &
$q^{-4}+q^{-3}+q^{-2}$ &$-q^{-2}-q^{-1}$ & 1 &0\\[0.5mm]
\normalsize $-$3& $q^{-6}+q^{-5}+2q^{-4}+q^{-3}+q^{-2}$ &
$-q^{-3}-q^{-2}-q^{-1}$ & 1  &0&0\\[0.5mm]
\normalsize $-$4& $-q^{-4}-q^{-3}-q^{-2}-q^{-1}$ & 1&0&0&0\\[0.5mm]
\normalsize $-$5& 1 &0&0&0&0}
Just as we did with the binomial coefficients, we must now interpret
our generalized Gaussian coefficients according to their combinatorial
propreties. We will do so in terms of
{\em linear partitions}. 
One usually defines a (linear) partition or a (linear) partition with
distinct parts to be a finite
sequence $\lambda $ of nonnegative integers with $\lambda _{1}\geq
\lambda _{2}\geq \cdots \geq \lambda _{k}$ or $\lambda _{1} > 
\lambda _{2} > \ldots > \lambda _{k} $ respectively. Since we are
now giving negative integers a status equal to that of
nonnegative ones, the classical definition becomes inappropriate.
We now procede to generalize the notion of a partition by first
generalizing the relations $\geq $ and $>$.

\begin{defn}[Hybrid Inequalities]
Given integers $i$ and $j$. We say that $i\ll j$ or equivalently $j\gg
i$ if and only if $i\in \{0,1,..,j-1 \}$. Similarly, we say that
$i\ull j$  or equivalently $j \ugg i$ if and only if $j\in \{1,2,..,j
\}$. 
\end{defn}

Note that these relations have their usual sense for positive
integers,  but behave rather oddly for negative integers.
For example, $-2,-1,0$ are the integers $\ull -3$, and
$-3,-2,-1$ are the integers $\ll -3$. 
These relations are not partial orders since they do not obey
reflexivity. However, they do obey antisymmetry and transitivity, so
we can still discuss their Hasse diagrams. 
\begin{figure}[htbp]
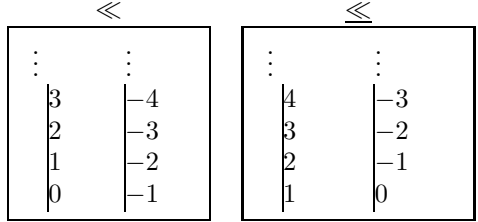

\label{fig:Hasse}
\caption{Hasse diagrams of $\ll$ and $\ull$}
\begin{center}
\begin{tabular}{cc}
$\ll $ &$\ull $\\
\fbox{\begin{tabular}{c|@{}p{8mm}|@{}p{8mm}}
	\multicolumn{2}{l}{\vdots}	&\multicolumn{1}{@{}l}{\vdots} \\
&	3 & $-4$ \\
&	2 & $-3$ \\
&	1 & $-2$ \\
&	0 & $-1$ 
	\end{tabular}}	&
\fbox{\begin{tabular}{c|@{}p{1cm}|@{}p{1cm}}
	\multicolumn{2}{l}{\vdots}	&\multicolumn{1}{@{}l}{\vdots} \\
&	4 & $-3$ \\
&	3 & $-2$ \\
&	2 & $-1$ \\
&	1 & $0$ 
	\end{tabular}}	
\end{tabular}
\end{center}
\end{figure}
As you can see in figure \ref{fig:Hasse}, the negative integers are
now unrelated to the 
positive integers,\footnote{$\ll $ and $\ull$ disagree on whether zero
is positive or negative.} and the usual order of negative integers has
been reversed. 

\begin{defn}[New Partitions]
Given a sequence $\lambda $ of integers, we say that $\lambda $  is an
{\em n-partition} of length $k$ (and width $n$) if $(n\ugg)\lambda
_{1}\ugg \lambda _{2} \ugg \cdots \ugg \lambda _{k}   $. We say 
that $\lambda $ is a {\em d-partition} of length $k$ (and width n) if
$(n \gg) \lambda_{1} \gg \lambda _{2} \gg \cdots \gg 
\lambda _{k}$. 

If $t$ is the sum of the $k$ parts of $\lambda $ then
we write $\lambda \vdash t$ or $|\lambda | = t$, and we say $\lambda $
is an n-partition or d-partition of $t$.  The length of $\lambda $ is
denoted $\ell(\lambda )$.
\end{defn}

By considering the case $t$ positive, we see that the n-partition generalizes
the idea of a {\bf n}ormal partition and the d-partition generalizes
the idea of a partition with {\bf d}istinct parts.\footnote{Note
however that for $t$ negative, it is the n-partitions which have
distinct parts instead of the d-partitions.}

For example, $(5,3,2,0)$ is a d-partition of 10 with length 4 and
minimum width 6. $(5,5,3,2)$ is an n-partition of 15 with length 4 and
minimum width 5.  $(-32,-7,-7)$ is a d-partition of $-46$ with length
3 and minimum (maximum?) width $-32$. $(-12,-7,-6,-3,0)$ is an
n-partition of $-22$ with length 5 and minimum (maximum?) width $-13$. 
Note also that a sequence can be at the same time a d-partition and an
n-partition, for example $(5,3,2)$.

It is well known (see for example \cite{Subspaces}) that for $n$ and
$k$ nonnegative, ${n\choose k}_{q}$ is a monic polynomial. Its
coefficient of $q^{t}$ is the number of partitions $\lambda$
of $t$ of length $n-k$ (or equivalently $k$) with all parts less than
or equal to $n$. 

In general, we now have:
\begin{thm}
Let $k$ be a nonnegative integer, and $n$ and $t$ be
arbitrary integers. Then ${n\choose k}_{q}=\sum _{t} c_{t}
q^{t}$ where  $c_{t}$ is the number of d-partitions $\lambda$
of $t$ with length $k$, and width $n$.
\end{thm}

{\em Proof:} By equation \ref{heq}.$\Box $

We can now express the $q$-binomial coefficients of regions 1 and 3 in
terms of each other.

\begin{cor}\label{duck}
Let $n$ and $k$ be nonnegative integers, then
$$ {n\choose k}_{q} = (-1)^{n+1}
q^{{n+1\choose 2 }-{k+1\choose 2 }} {-k-1 \choose -n-1}_{q}. \Box $$  
\end{cor}

\section{Stirling Numbers} \label{sns}
\subsection{Classical Region}

Our next example of connection constants is that of Stirling number of
the first
\fig{Stirling Numbers of the First Kind, $s(n,k)$---Region 1}{{c|rrrrrrr}
\label{fig:s1r1}
\normalsize $n \backslash k$& 
\normalsize 0 & \normalsize 1 & \normalsize 2 &
\normalsize 3 & \normalsize 4 & \normalsize 5 & \normalsize 6 \\
\hline
\normalsize 6& 0&$-120$& 274 & $-225$ & 85 & $-15$ & 1\\[0.5mm]
\normalsize 5& 0& 24 & $-50$ & 35 & $-10$ & 1 & 0\\[0.5mm]
\normalsize 4& 0& $-6$ & 11 & $-6$ & 1 & 0&0\\[0.5mm]
\normalsize 3& 0& 2 & $-3$ & 1 & 0&0&0\\[0.5mm]
\normalsize 2& 0& $-1$ & 1 & 0&0&0&0\\[0.5mm]
\normalsize 1& 0 & $1$ & 0&0&0&0&0\\[0.5mm]
\normalsize 0 & 1 & 0&0&0&0&0&0}
 and second kind,  $s(n,k)$ and $S(n,k)$ (see Figures \ref{fig:s1r1}
and \ref{s2r1}), 
\fig{Stirling Numbers of the Second Kind, $S(n,k)$---Region 1}{{c|rrrrrrrrrr}
\label{s2r1}
\normalsize $n \backslash k$& 
\normalsize 0 & \normalsize 1 & \normalsize 2 &
\normalsize 3 & \normalsize 4 & \normalsize 5 & \normalsize 6\\
\hline
\normalsize 6& 0& 1&31&90&65&15&1\\[0.5mm]
\normalsize 5& 0& 1&15&25&10&1 &0\\[0.5mm]
\normalsize 4& 0& 1&7&6&1& 0&0\\[0.5mm]
\normalsize 3& 0& 1&3&1 &0&0&0\\[0.5mm]
\normalsize 2& 0& 1&1 & 0&0&0&0\\[0.5mm]
\normalsize 1& 0& 1 & 0&0&0&0&0\\[0.5mm]
\normalsize 0& 1 & 0&0&0&0&0&0}
relaying the 
nonnegative powers of $x$  $(x^{n})_{n\geq 0}$ and the lower factorial
polynomials $((x)_{n})_{n\geq 0} = (x(x-1)\cdots (x-n+1))_{n\geq 0}$
\begin{eqnarray*}
(x)_{n} &=& \sum _{k=0}^{n}s(n,k) x^{k},\\
x^{n}&=& \sum _{k=0}^{n} S(n,k) (x)_{k}.
\end{eqnarray*}
Their best known combinatorial interpretations is that $S(n,k)$ counts
the number of partitions of an $n$-element set into $k$ blocks while
(up to sign) $s(n,k)$ counts 
the number of partitions of an $n$-element set containing $k$-cycles.

When generalizing these constants to negative values of $n$, we 
define the lower factorial sequence of rational functions
$((x)_{n})_{n\in \SB{Z}}$
for $n$ an arbitrary integer as the sequence of rational functions
with persistant roots $(i-1)_{i\in \SB{Z}}$. Next, we must
decide as we did in Proposition \ref{pr3} whether to consider
expansion of $(x)_{n}$ as a formal power series or as an inverse
Laurent series. In Proposition \ref{pr3}, the difference was not
critical; the binomial coefficients in regions 2 and 3 are related by
Proposition \ref{compl}. Each region is the reflection of the other.

\subsection{Formal Power Series Region}

Here, as we will see, there is no such symmetry. First, let us 
examine regions 1 and 2 defined by formal power series (see table
\ref{harmfig}). 
\fig{Stirling Numbers of the First Kind $s(n,k)$---Region 2}{{c|cccccr} 
\label{harmfig}
\normalsize $k \backslash n$& \normalsize $-5$ &
\normalsize $-4$ & \normalsize $-3$ & \normalsize $-2$ & \normalsize $-1$\\
\hline
\normalsize 6& $\frac{226,576,031,859}{5,598,720,000,000}$&
$\frac{11,679,655}{71,663,616}$ & $\frac{137,845}{279,936}$ &
$\frac{127}{128}$&1\\[2mm] 
\normalsize 5& $-\frac{3,673,451,957}{93,312,000,000}$ &
$-\frac{952,525}{5,971,968}$ & $-\frac{22,631}{46,656}$ & $-\frac{63}{64}$ &
$-1$ \\[2mm]
\normalsize 4& $\frac{58,067,611}{1,555,200,000}$ & $\frac{76,111}{497,664}$ &
$\frac{3661}{7776}$ & $\frac{31}{32}$ & 1 \\[2mm]
\normalsize 3& $-\frac{874,853}{25,920,000}$ & $-\frac{5845}{41,472}$ &
$-\frac{575}{1296}$ & $-\frac{15}{16}$ & $-1$ \\[2mm]
\normalsize 2& $\frac{12,019}{432,000}$ & $\frac{415}{3456}$ & $\frac{85}{216}$
& $\frac{7}{8}$ & 1 \\[2mm]
\normalsize 1& $-\frac{137}{7200}$ & $-\frac{25}{288}$ & $-\frac{11}{36}$ &
$-\frac{3}{4}$ & $-1$ \\[2mm]
\normalsize 0& $\frac{1}{120}$ & $\frac{1}{24}$ & $\frac{1}{6}$ & $\frac{1}{2}$
& 1 }
These
constants were considered in \cite{LR,ch2}. However, except for the
classical region 1, we no longer have integers, so there is no
combinatorial interpretation to be expected.
However, we do have several interesting recursions and identities
which prove to be useful in \cite{ILA} such as
\begin{thm} \cite[Porism 4.2 and Proposition 4.12]{ch2}
For all integers $n$ (positive or negative) and nonnegative integers $k$, 
$$ s(n,k)=\lim_{\epsilon \rightarrow 0}\Gamma (1-n+\epsilon)^{-1}
\comp _{k}\left(\frac{1}{\epsilon },\frac{1}{-1+\epsilon 
},\frac{1}{-2+\epsilon }, .. ,\frac{1}{-n+\epsilon} \right) .\Box $$
\end{thm}
We also have the recursion 
\begin{thm} \cite[Theorem 3.2]{ch2}
For all integers $n$ and for all positive integers $k$,
$$ s(n+1,k)=s(n,k-1)-ns(n,k).\Box  $$ 
\end{thm}
In region 2, we have the following result due to Knuth.
\begin{prop}\label{helpful} \cite[Proposition 4.3]{ch2}
Let $n$ and $k$ be nonnegative integers (not
both zero). Then $s(n,k)$ is given by the 
following finite sum:
\begin{equation}\label{HelpEq}
s(-n,k)= \frac{ (-1)^{k+1}}{n!} \sum_{m=1}^{n} {n\choose m} (-1)^{m} m^{-k}.
\end{equation}
\end{prop}

Note that Proposition~\ref{helpful} is the analog of the following classical
result \cite{Char} involving Stirling numbers of the {\em second} kind
$S(k,n)$.
$$ S(k,n)=\frac{(-1)^{n}}{n!}\sum_{m=1}^{k} {n\choose m}(-1)^{m}m^{k}. $$
Thus, in some sense to be made clear later, we can say that 
$$S(k,n)=(-1)^{n+k+1}s(-n,-k).$$

\subsection{Inverse Laurent Series Region}

On the other hand, in the case of inverse Laurent series, we can apply
the results of section \ref{TCSF} to calculate $s(n,k)$
\fig{Stirling Numbers of the First Kind, $s(n,k)$---Region 3}{{c|rrrrrr}
\label{fig:s1r3}
\normalsize $n \backslash k$& 
\normalsize $-6$ & \normalsize $-5$ & \normalsize $-4$ &
\normalsize $-3$ & \normalsize $-2$ & \normalsize $-1$ \\[0.5mm]
\hline
\normalsize $-1$& $-1$ & 1 & $-1$ & 1 & $-1$ & 1 \\[0.5mm]
\normalsize $-2$& 31 & $-15$ & 7 & $-3$ & 1 &0\\[0.5mm]
\normalsize $-3$& $-90$ & 25 & $-6$ & 1 & 0&0\\[0.5mm]
\normalsize $-4$& 65 & $-10$ & 1 & 0&0&0\\[0.5mm]
\normalsize $-5$& $-15$ & 1 &0&0&0&0\\[0.5mm]
\normalsize $-6$& 1 & 0&0&0&0&0}
	and $S(n,k)$. See tables \ref{fig:s1r3} and \ref{fig:s2r3}.
\fig{Stirling Numbers of the Second Kind, $S(n,k)$---Region 3}{{c|rrrrrr}
\label{fig:s2r3}\normalsize $n \backslash k$& 
\normalsize $-6$ & \normalsize $-5$ & \normalsize $-4$ &
\normalsize $-3$ & \normalsize $-2$ & \normalsize $-1$ \\
\hline
\normalsize $-1$& 120&24&6&2&1&1\\[0.5mm]
\normalsize $-2$& 274&50&11&3&1&0\\[0.5mm]
\normalsize $-3$& 225&35&6&1&0&0\\[0.5mm]
\normalsize $-4$& 85&10&1&0&0&0\\[0.5mm]
\normalsize $-5$& 15&1&0&0&0&0\\[0.5mm]
\normalsize $-6$& 1&0&0&0&0&0}
\begin{eqnarray*}
s(n,k) &=& \comp _{n-k} (0..n-1)\\
&=& \left\{ \begin{array}{ll}
(-1)^{n-k}e_{n-k}(1,2,\ldots,n-1)&\mbox{for $n\geq k\geq 0$,}\\[0.1in]
h_{n-k}(-1,-2,\ldots,n)		&\mbox{for $0>n\geq k$.}
         \end{array}
\right. \\[0.1in]
S(n,k) &=& \comp _{n-k} (-\{ 0..k\})\\
&=& \left\{ \begin{array}{ll}
h_{n-k}(1,2,\ldots,k)		&\mbox{for $n\geq k\geq 0$,}\\[0.1in]
(-1)^{k}e_{n-k}(-1,-2,\ldots,k+1)&\mbox{for $0>n\geq k$.}
         \end{array}
\right.
\end{eqnarray*}
By Theorem \ref{t5}, the numbers in regions 1 and 3 obey the identities 
\begin{eqnarray*}
s(n,k) &=& s(n-1,k-1) - n s(n-1,k)\\
S(n,k) &=& S(n-1,k-1) + k S(n-1,k).
\end{eqnarray*}

These Stirling numbers are always integers, so we can hope to give an
extended 
combinatorial interpretation. We will not only do that, but in fact,
we will give a combinatorial interpretation of the extended
$p,q$-Stirling numbers. The $p,q$-Stirling numbers of the first
\fig{$p,q$-Stirling Numbers of the First Kind, $s_{pq}(n,k)$---Region
1}{{r|*{5}{l}}    
\label{pq1r1}
\normalsize $n \backslash k$ & \normalsize 0 &\normalsize  1
&\normalsize  2 &\normalsize  3 &\normalsize  4\\ 
\hline  
\normalsize 4 & 0 & 
$-\sqrt{pq}^{3} [ \begin{array}[t]{@{}c@{}l@{}}
	q^{3}+2q^{2}p\\ +2qp^{2}-p^{3} &] \end{array}$
& $pq [ \begin{array}[t]{@{}c@{}l@{}}
      	(q+q^{2}+q^{3}) \\ +(1+q+2q^{2})p \\ +(1+2q)p^{2}+p^{3} &] 
        \end{array}$
& $-\sqrt{pq}(1+q+q^{2})$ & 1 \\[1mm]
\normalsize 3 & 0 & $pq[q+p]$& $-\sqrt{pq}(1+p+q)$ & 1 &0\\[0.5mm]
\normalsize 2 & 0 & $-\sqrt{pq}$ & 1 &0&0\\[0.5mm]
\normalsize 1 & 0 & 1 & 0&0&0\\[0.5mm]
\normalsize 0 & 1 & 0&0&0&0}
\fig{$p,q$-Stirling Numbers of the First Kind, $s_{pq}(n,k)$---Region
3}{{r|llll}   
\label{pq1r3}
\normalsize $n \backslash k$ &\normalsize  $-$4
&\normalsize  $-$3 &\normalsize  $-$2 &\normalsize  $-$1 \\
\hline  
\normalsize $-$1& $-1/\sqrt{pq}^{3}$ & $1/pq$ &
$-1/\sqrt{pq}$ & 1\\[0.5mm] 
\normalsize $-$2&
$[ \begin{array}[t]{@{}c@{}l@{}}
     p^{-2}+(q^{-1}+2)p^{-1}\\  +(1+q^{-1}+q^{-2}) & ] / pq
     \end{array}$ 
&$-[ \begin{array}[t]{@{}c@{}l@{}}
      p^{-1} \\  +(q^{-1}+1) & ]/ \sqrt{pq}
     \end{array} $ 
& 1 &0\\[1mm]
\normalsize $-$3&
$-[ \begin{array}[t]{@{}c@{}l@{}}
	p^{-2}+(q^{-1}+1)p^{-1} \\
	+(q^{-2}+q^{-1}+1) & ]/ \sqrt{pq}
	\end{array}$ 
& 1 &0&0\\[1mm]
\normalsize $-$4& 1&0&0&0}
kind (see tables \ref{pq1r1} and \ref{pq1r3}) and second kind
$s_{pq}(n,k)$ and  $S_{pq}(n,k)$ (see tables \ref{pq2r1} and
\ref{pq2r3})  are the 
coefficients connecting the powers of $x$ with the $p,q$-lower factorial
$(x;p,q)_{n}$ where the $p,q$-{\em lower factorials} $((x;p,q)_{n})_{n\in
\SB{Z}}$  form the sequence of rational functions with persistant
roots given by the $p,q$-{\em numbers}
$(\sqrt{pq}(p^{i}-q^{i})/(p-q))_{i\in \SB{Z}}$.\footnote{The
$p,q$-numbers are sometimes defined as $(p^{i}-q^{i})/(p-q)$ or simply
$p^{i}-q^{i}$. These alternate definitions lead to quite similar
results.} 
\fig{$p,q$-Stirling Numbers of the Second Kind, $S_{pq}(n,k)$---Region
1}{{r|*{6}{l}}   
\label{pq2r1}
\normalsize $n \backslash k$ & \normalsize
0 &\normalsize  1 &\normalsize  2 &\normalsize  3 &\normalsize  4\\
\hline  
\normalsize 4 & 0 & $\sqrt{pq}^{3}$ & 
$pq[ \begin{array}[t]{@{}c@{}l@{}}
     (1+q+q^{2}) \\ +(1+2q)p+p^{2} & ] \end{array}$ 
& $\sqrt{pq}[ \begin{array}[t]{@{}c@{}l@{}}
     (1+q+q^{2}) \\ +(1+q)p+p^{2} & ] \end{array}$ 
 & 1 \\[1mm]
\normalsize 3 & 0 & $pq$ & $\sqrt{pq}(1+q+p)$ & 1 &0\\[0.5mm]
\normalsize 2 & 0 & $\sqrt{pq}$ & 1 &0&0\\[0.5mm]
\normalsize 1 & 0 & 1 & 0&0&0\\[0.5mm]
\normalsize 0 & 1 & 0&0&0&0}
\fig{$p,q$-Stirling Numbers of the Second Kind, $S_{pq}(n,k)$---Region
3}{{r|llll}   
\label{pq2r3}
\normalsize $n \backslash k$ &\normalsize  $-$4
&\normalsize  $-$3 &\normalsize  $-$2 &\normalsize  $-$1 \\
\hline  
\normalsize $-$1&
$[ p^{-3}+2q^{-1}p^{-2} + 2q^{-2}p^{-1} +q^{-3}]/\sqrt{pq}^{3}$&
$[p^{-1}+q^{-1}]/pq$ & $1/\sqrt{pq}$ & 1\\[0.5mm] 
\normalsize $-$2&
$[ \begin{array}[t]{@{}c@{}l@{}}
     p^{-3}+(2q^{-1}+1)p^{-2} \\
     +(2q^{-2}+q^{-1}+1)p^{-1}\\
     + (q^{-3}+q^{-2}+q^{-1})& ]/ pq \end{array}$ 
&$ [ \begin{array}[t]{@{}c@{}l@{}}
     p^{-1} \\+(q^{-1}+1)& ]/\sqrt{pq} \end{array}$ 
& 1 &0\\[1mm]
\normalsize $-$3&
$[ \begin{array}[t]{@{}c@{}l@{}}
    p^{-2}+(q^{-1}+1)p^{-1}\\ +(q^{-2}+q^{-1}+1) & ] / \sqrt{pq} \end{array}$
& 1 &0&0\\[1mm]
\normalsize $-$4& 1&0&0&0}

The $p,q$-Stirling number obey the following properties:
\begin{eqnarray}
s_{pq}(n,k) &=& \comp _{n-k} ([0]_{pq}..[n-1]_{pq}) \nonumber \\
&=& \left\{ \begin{array}{ll}
(-1)^{n-k}e_{n-k}([1]_{pq},[2]_{pq},\ldots,[n-1]_{pq})& \mbox{for
$n\geq k\geq 0$,}\\[0.1in] 
h_{n-k}([-1]_{pq},[-2]_{pq},\ldots,[n]_{pq})		&\mbox{for
$0>n\geq k$.} 
         \end{array}
\right. \label{hepq}\\[0.1in]
S_{pq}(n,k) &=& \comp _{n-k} (-\{ [0]_{pq}..[k]_{pq}\}) \nonumber \\
&=& \left\{ \begin{array}{ll}
h_{n-k}([1]_{pq},[2]_{pq},\ldots,[k]_{pq})		&\mbox{for
$n\geq k\geq 0$,}\\[0.1in] 
(-1)^{k}e_{n-k}([-1]_{pq},[-2]_{pq},\ldots,[k+1]_{pq})&\mbox{for $0>n\geq k$.}
         \end{array}
\right. \label{HEPQ} \\
s_{pq}(n,k) &=& s_{pq}(n-1,k-1) - [n]_{pq} s_{pq}(n-1,k) \nonumber \\
S_{pq}(n,k) &=& S_{pq}(n-1,k-1) + [k]_{pq} S_{pq}(n-1,k). \nonumber 
\end{eqnarray}

By definition, the $p,q$-Stirling numbers transform into the
classical Stirling numbers by setting $p$ and $q$ equal to one. 
If we merely set $p$ equal to one, we have the usual $q$-Stirling
numbers $s_{q}(n,k)$
\fig{$q$-Stirling Numbers of the First Kind, $s_{q}(n,k)$---Region
1}{{r|*{5}{l}}    
\label{q1r1}
\normalsize $n \backslash k$ & \normalsize 0 &\normalsize  1
&\normalsize  2 &\normalsize  3 &\normalsize  4\\ 
\hline  
\normalsize 5 & 0 & 
$ \begin{array}[t]{@{}c@{}}
     q^{2}+3q^{3}+5q^{4}+6q^{5} \\ +5q^{6}+3q^{7}+q^{8}\end{array}$
 & $-[ \begin{array}[t]{@{}c@{}l@{}}
     4q^{3/2}+9q^{5/2}+12q^{7/2}+12q^{9/2} \\
     8q^{11/2}+4q^{13/2}+q^{15/2} & ] \end{array}$ 
 & $ \begin{array}{@{}c@{}}
     6q+9q^{2}+9q^{3} \\ +7q^{4}+3q^{5}+q^{6} \end{array}$
 & $-[ \begin{array}[t]{@{}c@{}l@{}}
     4q^{1/2}-3q^{3/2} \\+2q^{5/2}+q^{7/2} & ] \end{array}$\\[1mm]
\normalsize 4 & 0 & $-[q^{3/2}+2q^{5/2} +2q^{7/2}+q^{9/2}]$ &
$3q+4q^{2}+3q^{3}+q^{4}$ & $-[3q^{1/2}+2q^{3/2}+q^{5/2}]$ & 1 \\[0.5mm]
\normalsize 3 & 0 & $q+q^{2}$& $-[2q^{1/2}+q^{3/2}]$ & 1 &0\\[0.5mm]
\normalsize 2 & 0 & $-q^{1/2}$ & 1 &0&0\\[0.5mm]
\normalsize 1 & 0 & 1 & 0&0&0\\[0.5mm]
\normalsize 0 & 1 & 0&0&0&0}
\fig{$q$-Stirling Numbers of the First Kind, $s_{q}(n,k)$---Region
3}{{r|lllll}   
\label{q1r3}
\normalsize $n \backslash k$ &\normalsize $-5$  &\normalsize  $-$4
&\normalsize  $-$3 &\normalsize  $-$2 &\normalsize  $-$1 \\
\hline  
\normalsize $-$1& $q^{-2}$ & $-\sqrt{q}^{-5}$ & $q^{-1}$ &
$-\sqrt{q}^{-1}$ & 1\\[0.5mm] 
\normalsize $-$2& $- [ \begin{array}[t]{@{}c@{}l@{}}
      q^{-9/2}+4q^{-7/2}\\ +6q^{-5/2}+4q^{-3/2} & ] \end{array}$ 
& $\begin{array}[t]{@{}c@{}} q^{-3}+3q^{-2} \\ +3q^{-1} \end{array}$
&$-[ \begin{array}[t]{@{}c@{}l@{}}  q^{-3/2} \\-2q^{-1/2} & ] \end{array} $
& 1 &0\\[1mm]
\normalsize $-$3& 
$\begin{array}[t]{@{}c@{}}
    q^{-5}+3q^{-4}+7q^{-3} \\ +8q^{-2}+6q^{-1} \end{array}$ 
& $-[\begin{array}[t]{@{}c@{}l@{}}4
    q^{-5/2}+2q^{-3/2}\\+3q^{-1/2}&]\end{array}$ 
& 1&0&0\\[1mm]
\normalsize $-$4& $-[\begin{array}[t]{ll}
    q^{-7/2}-2q^{-5/2} \\ -3q^{-3/2}-4q^{-1/2} \end{array}$
 & 1&0&0&0\\[1mm]
\normalsize $-$5& 1 &0&0&0&0}
(see tables \ref{q1r1} and \ref{q1r3}) and $S_{q}(n,k)$ (see tables
\ref{q2r1} and \ref{q2r3}) defined as the coefficients connecting the
powers of $x$ with the $q$-lower factorial $(x;q)_{n}$ where the
$q$-lower 
factorials $((x;q)_{n})_{n\in \SB{Z}}$ form the sequence of rational
functions with persistant roots given by the $q$-numbers
$\sqrt{q}(q^{i}-1)/(q-1)$. 
\fig{$q$-Stirling Numbers of the Second Kind, $S_{q}(n,k)$---Region
1}{{r|*{6}{l}}   
\label{q2r1}
\normalsize $n \backslash k$ & \normalsize
0 &\normalsize  1 &\normalsize  2 &\normalsize  3 &\normalsize  4\\
\hline  
\normalsize 5 & 0 & $q^{2}$ 
& $\begin{array}[t]{@{}c@{}} 4q^{3/2}+6q^{5/2}\\ +4q^{7/2}+q^{9/2} \end{array}$ & 
$\begin{array}[t]{@{}c@{}} 6q+8q^{2}+7q^{3}\\+3q^{4}+q^{5}\end{array}$  &
$\begin{array}[t]{@{}c@{}}4q^{1/2}+3q^{3/2} \\+2q^{5/2}+q^{7/2}\end{array}$\\[1mm]
\normalsize 4 & 0 & $q^{3/2}$ & $3q+3q^{2}+q^{3}$ &
$3q^{1/2}+2q^{3/2}+q^{5/2}$ & 1 \\[0.5mm] 
\normalsize 3 & 0 & $q$ & $2q^{1/2}+q^{3/2}$ & 1 &0\\[0.5mm]
\normalsize 2 & 0 & $q^{1/2}$ & 1 &0&0\\[0.5mm]
\normalsize 1 & 0 & 1 & 0&0&0\\[0.5mm]
\normalsize 0 & 1 & 0&0&0&0}
\fig{$q$-Stirling Numbers of the Second Kind, $S_{q}(n,k)$---Region
3}{{r|lllll}   
\label{q2r3}
\normalsize $n \backslash k$ &\normalsize $-5$  &\normalsize  $-$4
&\normalsize  $-$3 &\normalsize  $-$2 &\normalsize  $-$1 \\
\hline  
\normalsize $-$1&
$\begin{array}[t]{@{}c@{}}
   q^{-8}+3q^{-7}+5q^{-6}+6q^{-5}\\+5q^{-4}+3q^{-3}+q^{-2}\end{array}$ 
& $\begin{array}[t]{@{}c@{}}
   q^{-9/2}+2q^{-7/2}\\+2q^{-5/2}+q^{-3/2}\end{array}$
& $q^{-2}+q^{-1}$ & $q^{-1/2}$ & 1\\[1mm]
\normalsize $-$2&
$\begin{array}[t]{@{}c@{}}
   q^{-15/2}+4q^{-13/2}+8q^{-11/2}+12q^{-9/2}\\
   +12q^{-7/2}+9q^{-5/2}+4q^{-3/2}\end{array}$ & 
$\begin{array}[t]{@{}c@{}}
   q^{-4}+3q^{-3}\\+4q^{-2}+3q^{-1}\end{array}$ 
&$\begin{array}[t]{@{}c@{}} q^{-3/2}\\+q^{-1/2}\end{array}$ & 1 &0\\[1mm]
\normalsize $-$3& 
$\begin{array}[t]{@{}c@{}}
	q^{-6}+3q^{-5}+7q^{-4}\\+9q^{-3}+9q^{-2}+6q^{-1}\end{array}$ &
$\begin{array}[t]{@{}c@{}}
       q^{-5/2}+2q^{-3/2}\\ +3q^{-1/2} \end{array}$ 
& 1  &0&0\\[1mm]
\normalsize $-$4& $q^{-7/2}+2q^{-5/2}+3q^{-3/2}+4q^{-1/2}$ & 1&0&0&0\\
\normalsize $-$5& 1 &0&0&0&0}

Note also that the $p,q$-numbers are invariant
under the interchange of $p$ with $q$. Thus, the $p,q$-Stirling
numbers are symmetric in $p$ and $q$.

This interpretation will be made in terms of Ferrers diagrams.
\begin{defn}[Ferrers Diagram]
The {\em Ferrers diagram} of a n-partition or d-partition $\lambda $
with length $k$ 
is defined as the following hybrid set of squares $(i,j)$:
$$ {\cal H}(\lambda ) = \sum _{i=1}^{k} \{ (1,1) .. (1,\lambda _{i})
\}. $$ 
\end{defn}
The cardinality of ${\cal H}(\lambda )$ is $|\lambda |.$ 

In the case of classical partitions, the above corresponds to the
traditional concept of a Ferrers diagram. For example, the Ferrers
diagram of (3,2,0) is the positive set
$\{(1,1),(1,2),(1,3),(2,1),(2,2)| \}$. On the other hand, the Ferrers
diagram of $(-3,-2)$ is the negative set
$\{|(1,0),(1,-1),(1,-2),(2,0),(2,-1) \}$. 

\begin{defn}[0-1 Tableaux]
A {\em 0-1 p-tableau} (resp. a {\em 0-1 n-tableau}) is a pair
$\alpha= (\lambda ,f)$ where $\lambda $ is a p-partition (resp. n-partition),
and $f$ is a  {\em filling} of that diagram. That is to say it 
is a collection of constants equal to 0 or 1 and denoted $f_{ij}$ for
each $(i,j)\in {\cal H}(\lambda )$ subject to the condition that there is
exactly one 1 in each row. ie: for all $1\leq i\leq \ell (\lambda )$,
there is a unique $j$ such that $f_{ij}=1$.  

The {\em inversion number} $\inv(\alpha)$ of a 0-1 tableau
$\alpha=(\lambda ,f)$ is one half plus the number of zeros to 
the left of a one in the diagram, ie: the number of triples
$(i,j_{1},j_{2})$ such that $j_{1}\ll j_{2}$, $f_{ij_{1}}=0$, and
$f_{ij_{2}}=1$ while taking into account possible negative
multiplicities. Similarly, the {\em non-inversion number}
$\nin(\alpha)$ of a 0-1 tableau is one half plus the number of zeros
to the right of a one in the diagram. 
\end{defn}

We then have the following generalization of \cite[eq. (1.4)]{MR}:

\begin{thm}
Let $n$ and $k$ be integers. Then the Stirling numbers of the second
kind of degree $n$ and order $k$ are given by sums
\begin{eqnarray*}
 s_{pq}(n,k) &= &\sum _{\alpha} q^{\inv (\alpha)} p^{\nin (\alpha)}\\
 s_{q}(n,k) &= &\sum _{\alpha} q^{\inv (\alpha)}
\end{eqnarray*}
over 0-1 tableaux $\alpha$ of d-partitions of width $n$ and length
$n-k.$ Thus, $s(n,k)$ is the number of such tableaux.

The Stirling numbers of the first kind of degree $n$ and order $k$ are
given by similar sums 
\begin{eqnarray*}
 S_{pq}(n,k) &=& \sum _{\alpha} q^{\inv (\alpha)} p^{\nin (\alpha)}\\
 S_{q}(n,k) &=& \sum _{\alpha} q^{\inv (\alpha)} p^{\nin (\alpha)}
\end{eqnarray*}
this time over 0-1  tableaux $\alpha$ of n-partitions of width $k$,
and length $n-k .$ Thus, $S(n,k)$ is the number of such tableaux. 
\end{thm}

{\em Proof:} Equation \ref{hepq}.$\Box $

\begin{cor}\label{inversion}
Let $n$ and $k$ be integers. Then we have the following relations
\begin{eqnarray*}
S_{pq}(n,k) &=& (-1)^{n+k} s_{p^{-1}q^{-1}}(-k,-n)\\
S_{q}(n,k) &=&  (-1)^{n+k} s_{q^{-1}}(-k,-n)\\
S(n,k) &=&  (-1)^{n+k} s(-k,-n).\Box 
\end{eqnarray*}
\end{cor}

 \end{document}